\documentclass[a4paper,12pt]{article}

\usepackage{color}

\usepackage{makeidx}

\usepackage{latexsym}
\usepackage{amscd}
\usepackage{amsmath}
\usepackage{amssymb}

\usepackage{amsthm}
\usepackage{float}
\usepackage{graphicx}
\usepackage{psfrag}

\theoremstyle{plain}
\newtheorem{theorem}{Theorem}
\newtheorem{example}[theorem]{Example}
\newtheorem{corollary}[theorem]{Corollary}
\newtheorem{lemma}[theorem]{Lemma}

\newtheorem{proposition}[theorem]{Proposition}
\newtheorem{conjecture}[theorem]{Conjecture}

\theoremstyle{definition}
\newtheorem{definition}[theorem]{Definition}
\newtheorem{remark}[theorem]{Remark}

\newdimen\argwidth
\def\db[#1\db]{%
 \setbox0=\hbox{$#1$}\argwidth=\wd0
 \setbox0=\hbox{$\left[\box0\right]$}
  \advance\argwidth by -\wd0
 \left[\kern.3\argwidth\box0 \kern.3\argwidth\right]}

\newcommand{\bC}{\ensuremath{\mathbb{C}}}

\newcommand{\bP}{\ensuremath{\mathbb{P}}}

\newcommand{\bR}{\ensuremath{\mathbb{R}}}

\newcommand{\bZ}{\ensuremath{\mathbb{Z}}}

\newcommand{\scA}{\ensuremath{\mathcal{A}}}

\newcommand{\scE}{\ensuremath{\mathcal{E}}}

\newcommand{\scO}{\ensuremath{\mathcal{O}}}

\newcommand{\frakm}{\ensuremath{\mathfrak{m}}}

\newcommand{\Der}{\mathrm{Der}}

\newcommand{\Proj}{\operatorname{Proj}}

\title{Coxeter multiarrangements with quasi-constant multiplicities}
\author{
Takuro Abe
 \ \ \ 
Masahiko Yoshinaga
}
\date{\today}
\begin{document}

\maketitle

\begin{abstract}
We study structures of derivation modules of 
Coxeter multiarrangements with quasi-constant 
multiplicities 
by using the primitive derivation. 
As an application, we show that the characteristic 
polynomial of a Coxeter multiarrangement with quasi-constant 
multiplicity is combinatorially computable. 
\end{abstract}

\section{Introduction}
\label{sec:intro}

Let $V$ be an $\ell$-dimensional Euclidean space over $\bR$ 
with inner product $I:V\times V\rightarrow \bR$. 
Fix a coordinate 
$(x_1, \cdots, x_\ell)$ and put 
$S=S(V^*)\otimes_\bR\bC=\bC[x_1, \ldots,x_\ell]$. 
Let $W\subset O(V, I)$ be a finite irreducible reflection group with 
the Coxeter number $h$. 
It is proved by Chevalley in \cite{ch} that the invariant ring 
$S^W$ is a polynomial ring $S^W=\bC[P_1, \ldots, P_\ell]$ 
with $P_1, \ldots, P_\ell$ are homogeneous generators. 
Suppose that $\deg P_1\leq\cdots\leq\deg P_\ell$. Then it is 
known that $\deg P_1=2<\deg P_2\leq\cdots\leq\deg P_{\ell-1}<\deg P_\ell=h$. 
Let $\scA$ be the corresponding Coxeter arrangement, i.e., 
the collection of all reflecting hyperplanes of $W$. 
Fix a defining linear form $\alpha_H\in V^*$ for each 
hyperplane $H\in\scA$. 
Let $m:\scA\rightarrow\bZ_{\geq0}$ be a map, called 
a \textit{multiplicity} on $\scA$. Then the 
pair $(\scA, m)$ is called a \textit{Coxeter multiarrangement}. 
Let $\Der(S)$ denote the module of $\bC$-linear 
derivations of $S$. Define a graded $S$-module $D(\scA,m)$ by 
$$
D(\scA, m)=
\{
\delta\in\Der(S)\mid \delta\alpha_H\in(\alpha_H)^{m(H)}\ \mbox{for all}\ H\in\scA
\}. 
$$
We say a multiarrangement $(\scA,m)$ is 
\textit{free} if $D(\scA,m)$ is a free 
$S$-module. When $(\scA,m)$ is free, we can choose a homogeneous basis 
$\{\theta_1,\ldots,\theta_\ell\}$ for $D(\scA,m)$ and 
call the multiset $(\deg(\theta_1),\ldots,\deg(\theta_\ell))$ the 
\textit{exponents} of a free multiarrangement $(\scA,m)$ and denoted by 
$\exp(\scA,m)$, where the degree is the polynomial degree. 
The module $D(\scA,m)$ 
was first defined by Ziegler (\cite{zie-multi}) and 
deeply studied for Coxeter multiarrangements with constant multiplicity 
by \cite{st-double, ter-multi}. 
In particular, Terao proved that if $m$ is constant, 
then $(\scA, m)$ is free and the exponents are expressed 
by using exponents of the Coxeter group and the Coxeter 
number $h$ (\cite{ter-multi}). 
These facts 
played a crucial role in the proof of Edelman-Reiner 
conjecture (\cite{er-rhomb, yos-char}).

Another aspect of the above module is a relation 
with the Hodge filtration of $\Der(S^W)$ 
introduced by K. Saito in \cite{sai-lin, sai-unif}. 
It is proved in \cite{ter-hodge} that if $m$ is 
a constant multiplicity with $m=2k+1$, then 
the $S^W$-module $D(\scA, m)^W$ of all $W$-invariant vector fields is 
precisely equal to the $k$-th Hodge filtration of $\Der(S^W)$. 
Based on these results, 
a geometrically expressed $S$-basis of the module $D(\scA, m)$ 
for special kind of (not necessarily constant) multiplicities was 
constructed in \cite{yos-multi}. 
The purpose of this paper is to strengthen and generalize 
results in \cite{ter-multi, yos-multi} by developing the ``dual'' 
version of 
\cite{yos-multi}. Indeed, we handle the following 
``quasi-constant'' multiplicities. 
\begin{definition}
A multiplicity $\widetilde{m}:\scA\rightarrow\bZ_{\geq 0}$ 
is said to be \textit{quasi-constant} if 
$$
\max\{\widetilde{m}(H)\mid H\in\scA\}-
\min\{\widetilde{m}(H)\mid H\in\scA\}\leq 1. 
$$
\end{definition}
It is clear that for a given quasi-constant multiplicity 
$\widetilde{m}$, there exist an integer $k$ and a 
$\{0,1\}$-valued multiplicity $m:\scA\rightarrow \{0,1\}$ 
such that $\widetilde{m}$ is either $2k+m$ or $2k-m$. 
The above $k\in\bZ_{\geq 0}$ and $m$ are uniquely determined 
unless $\widetilde{m}$ is 
the constant multiplicity with odd value. 
Our main results are concerning structures of derivation modules 
for Coxeter arrangements with quasi-constant multiplicities. 

\begin{theorem}
Let $\scA$ be a Coxeter arrangement with the Coxeter number 
$h$ and $m:\scA\rightarrow \{0,1\}$ be a $\{0,1\}$-valued 
multiplicity. Then 
\begin{itemize}
\item[(1)] $D(\scA, 2k+m)\cong D(\scA, m)(-kh)$, 
\item[(2)] $D(\scA, 2k-m)\cong \Omega^1(\scA, m)(-kh)$, and 
\item[(3)] The modules 
$D(\scA, 2k+m)(kh)$ and $D(\scA, 2k-m)(kh)$ are dual $S$-modules to each other,  
\end{itemize}
where $M(n)$ denotes the degree shift by $n$ for a 
graded $S$-module $M$. 
\label{seconemain}
\end{theorem}

Theorem \ref{seconemain} generalizes \cite{ter-multi, yos-multi}  in the following 
three parts. In \cite{yos-multi}, the isomorphism 
$D(\scA, 2k+m)\cong D(\scA, m)(-kh)$ is proved for the case 
$(\scA, m)$ is free. In Theorem \ref{seconemain}, the 
assumption on the freeness is removed. Furthermore, considerations on 
$\Omega^1(\scA, m)$ instead of $D(\scA, m)$ enable us to 
treat multiarrangements of the type $(\scA, 2k-m)$ as well (2). 
The structure of the module $D(\scA, m)$ is not so much known when 
it is not free. Combining Theorem \ref{seconemain} (1) and (2), 
we have an interesting relation Theorem \ref{seconemain} (3), i.e., there exists a natural 
pairing between the modules $D(\scA, 2k+m)$ and $D(\scA, 2k-m)$. 
It may be simply said that a 
relation between multiplicities 
gives an algebraic relation between derivation modules. 


The organization of this paper is as follows. In \S\ref{sec:terao} 
we review Terao's result about the derivation modules of 
Coxeter arrangements with constant 
multiplicity in \cite{ter-multi} 
from the viewpoint of 
the differential modules. In \S\ref{sec:main} we prove 
Theorem \ref{seconemain} (2) and the rest in \S\ref{sec:concl}. 
In \S\ref{sec:char} we apply these results to compute 
characteristic polynomials for Coxeter multiarrangements 
with quasi-constant multiplicities. 

\section{An interpretation of Terao's basis} 
\label{sec:terao}

In this section, we recall the main result of 
\cite{ter-multi} and 
give an interpretation 
through the dual 
basis for $\Omega^1(\scA, m)$. 
Let us first recall the definition of $\Omega^1(\scA, m)$. 

\begin{definition}
Put $Q(\scA, m)=\prod_{H\in\scA}\alpha_H^{m(H)}$ and denote by 
$\Omega_V^1=S\otimes_\bC V^*=\bigoplus_{i=1}^\ell S\cdot dx_i$ 
the module of differentials. Define 
$$
\Omega^1(\scA, m)=
\left\{\left.
\omega\in\frac{1}{Q(\scA, m)}\Omega_V^1
\right|
\begin{array}{cc}
d\alpha_H\wedge\omega\mbox{ does not have poles}\\
\mbox{ along $H$, 
for any $H\in\scA$} 
\end{array}
\right\}. 
$$
\end{definition}

It is known that $\Omega^1(\scA,m)$ is the dual $S$-module 
of $D(\scA,m)$ and vice versa (\cite{sai-log}, \cite{zie-multi}). 
Next we define the affine connection $\nabla$. 
\begin{definition}
For a given rational vector field 
$\delta=\sum_{i=1}^\ell f_i\frac{\partial}{\partial x_i}$ 
and a rational differential $k$-form 
$\omega=\sum_{i_1, \ldots, i_k}g_{i_1, \ldots, i_k} dx_{i_1, \ldots, i_k}$, 
define $\nabla_\delta\omega$ by 
$$
\nabla_\delta\omega=
\sum_{i_1, \ldots, i_k}\delta(g_{i_1, \ldots, i_k}) dx_{i_1, \ldots, i_k}. 
$$
\end{definition}
The above $\nabla$ defines a connection. 
We collect some elementary properties of $\nabla$ which 
will be used later. 

\begin{proposition}\label{prop:elem}
For a rational vector field $\delta$, rational differential 
form $\omega$ and $f \in S$, 
$\nabla$ has the following properties. 
\begin{itemize}
\item $\nabla_\delta f=\delta(f)$. 
\item $\nabla_{f\delta}\omega=f\nabla_\delta\omega$. 
\item Leibniz rule: $\nabla_{\delta}(f\omega)=f\nabla_\delta\omega+(\delta f)\omega$. 
\item For any linear form $\alpha\in V^*$, 
$\nabla_\delta(d\alpha\wedge\omega)=d\alpha\wedge\nabla_\delta\omega$. 
\end{itemize}
\end{proposition}

Now we fix a generating system $P_1, \ldots, P_\ell$ of the invariant 
ring $S^W=\bC[P_1, \ldots, P_\ell]$ as in \S\ref{sec:intro}. 
Note that we may choose $P_1(x)=I(x,x)$. 
Then $\frac{\partial}{\partial P_i}$ ($i=1, \ldots, \ell$) can be 
considered as a rational vector field on $V$ with order one 
poles along $H\in\scA$. Especially, 
we denote $D=\frac{\partial}{\partial P_\ell}$ and call it 
the \textit{primitive derivation}. 
Since $\deg P_i<\deg P_\ell$ for $i\leq\ell-1$, 
the primitive derivation $D$ is uniquely determined 
up to nonzero constant multiple 
independent of the choice of the generators $P_1, \ldots, P_\ell$ 
(\cite{sai-lin, sai-unif}). 

For any constant multiplicity $m\in\bZ_{\geq 0}$, Terao showed the 
freeness of $\Omega^1(\scA, m)$ by constructing a basis. 
\begin{theorem}\label{thm:ter}
\cite[Theorem 1.1]{ter-multi}
\begin{itemize}
\item[(1)] 
If 
$m=2k$, then 
$$
\nabla_{\frac{\partial}{\partial x_1}}\nabla_D^k dP_1, 
\nabla_{\frac{\partial}{\partial x_2}}\nabla_D^k dP_1, 
\ldots, 
\nabla_{\frac{\partial}{\partial x_\ell}}\nabla_D^k dP_1
$$
forms a basis for $\Omega^1(\scA, 2k)$. \\
\item[(2)] 
If 
$m=2k+1$, then 
$$
\nabla_{\frac{\partial}{\partial P_1}}\nabla_D^k dP_1, 
\nabla_{\frac{\partial}{\partial P_2}}\nabla_D^k dP_1, 
\ldots, 
\nabla_{\frac{\partial}{\partial P_\ell}}\nabla_D^k dP_1
$$
forms a basis for $\Omega^1(\scA, 2k+1)$. 
\end{itemize}
\end{theorem}

Originally in \cite{ter-multi} a basis for $D(\scA, m)$ is 
constructed. The above expression is obtained just by switching 
to $\Omega^1(\scA, m)$ through $\nabla$.

\section{Main results}
\label{sec:main}

\begin{lemma}\label{lem:indep}
Let $\delta_1, \ldots, \delta_\ell$ be rational vector fields. 
Suppose that they are linearly independent over $S$. Then 
$$
\nabla_{\delta_1}\nabla_D^k dP_1, 
\nabla_{\delta_2}\nabla_D^k dP_1, 
\ldots, 
\nabla_{\delta_\ell}\nabla_D^k dP_1
$$
are linearly independent over $S$. 
\end{lemma}
\proof 
Put $\delta_i=\sum_{j=1}^\ell a_{ij}\partial_j$, where 
$\partial_j=\frac{\partial}{\partial x_j}$. Then linearly independence 
of $\{\delta_1,\ldots,\delta_\ell\}$ 
is equivalent to $\det (a_{ij})\neq 0$. 
Now the assertion is clear from Theorem \ref{thm:ter} (1) and 
$$
\nabla_{\delta_i}\nabla_D^k dP_1= 
\sum_{j=1}^\ell a_{ij}\nabla_{\partial_j}\nabla_D^k dP_1. 
$$
\qed

\begin{lemma}\label{lem:transv}
The pole order of $\nabla_D^k dP_1$ is exactly equal to 
$2k-1$. More precisely, $\nabla_D^k dP_1\in\frac{1}{Q(\scA, 2k-1)}\Omega^1_V$ 
and $\alpha_H^{2k-2}\nabla_D^k dP_1$ has a pole along $H$ for any $H\in\scA$. 
\end{lemma}
\proof 
First note that since 
$$
\nabla_D^k dP_1=
\nabla_{\frac{\partial}{\partial P_\ell}}\nabla_D^{k-1} dP_1, 
$$
Theorem \ref{thm:ter} implies that 
$\nabla_D^k dP_1\in\Omega^1(\scA, 2k-1)$. 
Hence $\nabla_D^k dP_1\in\frac{1}{Q(\scA, 2k-1)}\Omega^1_V$. 

Suppose that there exists $H\in\scA$ such that 
$\alpha_H^{2k-2}\nabla_D^k dP_1$ does not have poles along $H$. 
Let us define the characteristic multiplicity $m_H$ by 
$$
m_H(H')=
\left\{
\begin{array}{cc}
1&\mbox{ if } H'=H,\\
0&\mbox{ if } H'\neq H. 
\end{array}
\right.
$$
Then it is easily seen that $\nabla_D^k dP_1\in\Omega^1(\scA, 2k-1-m_H)$. 
Since $\nabla_{\frac{\partial}{\partial P_j}}$ increases the pole order 
at most two, we have 
$\nabla_{\frac{\partial}{\partial P_j}}\nabla_D^k dP_1\in
\Omega^1(\scA, 2k+1-m_H)$. 
However, this contradicts to Theorem \ref{thm:ter} (2), for 
$\Omega(\scA, 2k+1)\supsetneqq \Omega^1(\scA, 2k+1-m_H)$. 
\qed

\begin{remark}
Lemma \ref{lem:transv} is a dual counterpart to \cite[Lemma 4]{yos-multi}. 
This property is related to the regularity of eigenvectors of the Coxeter 
element, which is of crucial importance in \cite{sai-lin, sai-unif}. 
\end{remark}

Let $m:\scA\rightarrow\{0,1\}$ be a $\{0,1\}$-valued multiplicity. 
The primitive derivation and $\nabla$ enable us to compare 
$D(\scA, m)$ and $\Omega^1(\scA, 2k-m)$.

\begin{theorem}\label{thm:isom}
For $\delta\in D(\scA, m)$, 
$\Phi_k(\delta):=\nabla_\delta\nabla_D^k dP_1$ is contained in 
$\Omega^1(\scA, 2k-m)$. Furthermore, 
the map 
$$
\begin{array}{cccc}
\Phi_k:&D(\scA, m)(kh)&\longrightarrow&\Omega^1(\scA, 2k-m)\\
&&&\\
&\delta&\longmapsto&\nabla_\delta\nabla_D^k dP_1
\end{array}
$$
gives an $S$-isomorphism. 
\end{theorem}

\proof 
Since $\nabla_\delta$ increases pole order at most one, from Lemma \ref{lem:transv}, 
$\nabla_\delta\nabla_D^k dP_1\in\frac{1}{Q(\scA, 2k)}\Omega_V^1$. 
Let $H\in\scA$. Then $\alpha_H^{2k-1}\cdot\nabla_D^k dP_1$ has no 
poles along $H$. Thus 
$\nabla_\delta\left(\alpha_H^{2k-1}\cdot\nabla_D^k dP_1\right)$ also has 
no poles along $H$. Suppose $m(H)=1$, and put 
$\delta(\alpha_H)=\alpha_H g$. 
Then we have 
\begin{eqnarray}
\nabla_\delta\left(\alpha_H^{2k-1}\cdot\nabla_D^k dP_1\right)&=&
(2k-1)\alpha_H^{2k-2}\delta(\alpha_H)\nabla_D^k dP_1+
\alpha_H^{2k-1}\nabla_\delta\nabla_D^k dP_1 \nonumber
\\
&=&
(2k-1)\alpha_H^{2k-1}g\nabla_D^k dP_1+
\alpha_H^{2k-1}\nabla_\delta\nabla_D^k dP_1. \nonumber
\end{eqnarray}
Hence $\alpha_H^{2k-1}\nabla_\delta\nabla_D^k dP_1$ has no pole 
along $H$. This shows that 
$\nabla_\delta\nabla_D^k dP_1\in\frac{1}{Q(\scA, 2k-m)}\Omega_V^1$. 
Since $d\alpha_H\wedge\nabla_D^k dP_1$ has no poles along $H$, 
using Proposition \ref{prop:elem}, 
$\nabla_\delta(d\alpha_H\wedge\nabla_D^k dP_1)=
d\alpha_H\wedge\nabla_\delta\nabla_D^k dP_1$ also does not have 
poles along $H$. This means 
$\Phi_k(\delta)=\nabla_\delta\nabla_D^k dP_1\in\Omega^1(\scA, 2k-m)$. 

Next we prove the injectivity. Let $K$ be the field of all 
rational functions. Since $\Phi_k$ is $S$-homomorphic, it can be 
extended to a $K$-linear map 
$$
\widetilde{\Phi_k}:
D(\scA, m)\otimes_SK\longrightarrow 
\Omega^1(\scA, 2k-m)\otimes_SK. 
$$
Then 
$\widetilde{\Phi_k}$ is isomorphic due to Lemma \ref{lem:indep}. 
Hence the induced map $\Phi_k$ is obviously injective. 

Finally we prove the surjectivity. 
Let $\omega\in\Omega^1(\scA, 2k-m)$. Then clearly 
$\omega\in\Omega^1(\scA, 2k)$. Hence from Theorem \ref{thm:ter}, 
there exists $\delta\in D(\scA, 0)=\sum_i S\partial_i$ such that 
$\omega=\nabla_\delta\nabla_D^k dP_1$. 
If $m\equiv 0$, there is nothing to prove. Otherwise, 
choose a hyperplane 
$H\in\scA$ such that $m(H)=1$. Then 
$\nabla_\delta\left(\alpha_H^{2k-1}\cdot\nabla_D^k dP_1\right)=
(2k-1)\alpha_H^{2k-2}\delta(\alpha_H)\nabla_D^k dP_1+
\alpha_H^{2k-1}\omega$ does not have poles along $H$. 
Hence 
$\alpha_H^{2k-2}\delta(\alpha_H)\nabla_D^k dP_1$ does not have 
poles along $H$. From Lemma \ref{lem:transv}, $\delta(\alpha_H)$ has 
to be divisible by $\alpha_H$. This shows that $\delta\in D(\scA, m)$. 
\qed

\section{Conclusions} 
\label{sec:concl}

By using parallel arguments to \S\ref{sec:main} in the context of \cite{yos-multi}, 
we can prove the following result. 
The notation is the same as above. 

\begin{theorem}
\label{thm:der}
Let $m:\scA\rightarrow\{0,1\}$ be a $\{0,1\}$-valued multiplicity and 
$E=\sum x_i\partial_i$ be the Euler vector field. Then for $\delta\in D(\scA, m)$, 
$\Psi_k(\delta):=\nabla_\delta\nabla_D^{-k}E$ is contained in 
$D(\scA, 2k+m)$. Furthermore, the map 
$$
\begin{array}{cccc}
\Psi_k:&D(\scA, m)(-kh)&\longrightarrow&D(\scA, 2k+m)\\
&&&\\
&\delta&\longmapsto&\nabla_\delta\nabla_D^{-k}E
\end{array}
$$
gives an $S$-isomorphism. 
\end{theorem}

The action of $\nabla_D$ shifts degree by $-h$. This proves the 
following results. 

\begin{corollary}
For a $\{0,1\}$-valued multiplicity $m:\scA\rightarrow\{0,1\}$ and 
an integer $k>0$, 
the following conditions are equivalent. 
\begin{itemize}
\item $(\scA, m)$ is free with exponents $(e_1, \ldots, e_\ell)$. 
\item $(\scA, 2k+m)$ is free with exponents $(kh+e_1, \ldots, kh+e_\ell)$. 
\item $(\scA, 2k-m)$ is free with exponents $(kh-e_1, \ldots, kh-e_\ell)$. 
\end{itemize}
\label{chordal}
\end{corollary}

\begin{remark}
The first condition 
in Corollary \ref{chordal} is equivalent to 
say the subarrangement $m^{-1}(1)\subset\scA$ is free. 
For the Coxeter arrangement of type $A$, free subarrangements 
$(\scA, m)$ are completely classified in \cite{sta-super}. 
See also \cite{er-ab}. 
\end{remark}


Another conclusion is the following. 

\begin{theorem}
\label{thm:dual}
Let $(\scA, m)$ be a Coxeter arrangement 
with a $\{0,1\}$-valued multiplicity $m$ and $k>0$. Then 
$D(\scA, 2k+m)(kh)$ and $D(\scA, 2k-m)(kh)$ are dual $S$-module 
to each other. 
\end{theorem}

\proof
Combining Theorem \ref{thm:isom} and \ref{thm:der}, we have 
the following isomorphisms of graded $S$-modules. 
$$
D(\scA, 2k+m)(kh)\cong D(\scA, m)\cong \Omega^1(\scA, 2k-m)(-kh). 
$$
Since $\Omega^1(\scA, 2k-m)\cong D(\scA, 2k-m)^*$, 
we have $D(\scA, 2k+m)(kh)\cong \Omega^1(\scA, 2k-m)(-kh)\cong D(\scA, 2k-m)(kh)^*$. 
\qed

\section{Characteristic polynomials}
\label{sec:char}

In a recent paper \cite{atw-char}, the characteristic 
polynomial $\chi((\scA, m), t)\in\bZ[t]$ for a 
multiarrangement $(\scA, m)$ is defined. In this 
section, we apply results in the previous sections to 
compute the characteristic polynomials. 
Let us first recall the definition of the characteristic 
polynomial briefly. 

Let $(\scA, m)$ be a multiarrangement of rank $\ell$. 
Then the module $D^p(\scA, m)$ and 
$\Omega^p(\scA,m)$ are defined for $0\leq p\leq \ell$ 
(see Introduction of \cite{atw-char} and \cite{zie-multi}), 
and define functions 
\begin{eqnarray*}
\psi(\scA, m; t,q)=\sum_{p=0}^\ell H(D^p(\scA, m), q)(t(q-1)-1)^p, \\
\phi(\scA, m; t,q)=\sum_{p=0}^\ell H(\Omega^p(\scA, m), q)(t(1-q)-1)^p,  
\end{eqnarray*}
in $t$ and $q$, where $H(M, q)$ is the Hilbert series of a graded $S$-module $M$. 
In \cite{atw-char}, $\psi$ and $\phi$ are proved to be polynomials in 
$t$ and $q$ and $(-1)^\ell\psi(\scA, m; t, 1)=\phi(\scA, m; t, 1)$. 
The characteristic polynomial of $(\scA, m)$ is by definition 
$$
\chi((\scA, m), t)=(-1)^\ell\psi(\scA, m; t, 1)=\phi(\scA, m; t, 1). 
$$
Note that the above definition is a generalization 
of so called Solomon-Terao's formula (\cite{st-stf}), that 
is, $\chi((\scA, 1), t)$ is equal to the combinatorially 
defined characteristic polynomial $\chi(\scA, t)$ of $\scA$ 
(\cite{ot-arr}). 

In general the computation of the characteristic polynomial 
$\chi((\scA, m), t)$, 
especially the constant term, is difficult. One of the 
reasons is that $\chi((\scA, m), t)$ is not a combinatorial 
invariant. However, we can compute it 
combinatorially for Coxeter multiarrangements 
with quasi-constant multiplicities.

\begin{theorem}
\label{thm:shift}
Let $\scA$ be a Coxeter arrangement with the Coxeter number 
$h$, and $m:\scA\rightarrow\{0,1\}$ 
be a $\{0,1\}$-valued multiplicity as in the previous sections. 
Let $k\in\bZ_{>0}$. Then 
\begin{itemize}
\item[(1)] $\chi((\scA, 2k+m), t)=\chi((\scA, m), t-kh)$, and 
\item[(2)] $\chi((\scA, 2k-m), t)=(-1)^\ell\chi((\scA, m), kh-t)$. 
\end{itemize}
\end{theorem}
For the proof, we need the following lemmas.

\begin{lemma}
\label{lem:sat}
Let $\frakm=(x_1, \ldots, x_\ell)\subset S$ be the graded 
maximal ideal of $S$. 
Let $(\scA, m)$ be any multiarrangement. Then 
$\Omega^p(\scA, m)$ is saturated in the following sense, 
that is, if $\omega\in\frac{1}{Q(\scA, m)}\Omega^p_V$ satisfies 
$\frakm\cdot\omega\subset\Omega^p(\scA, m)$, then 
$\omega\in\Omega^p(\scA, m)$. 
Similarly, if $\delta\in\Der^p(S)$ satisfies $\frakm\cdot\delta\subset D^p(\scA, m)$, 
then $\delta\in D^p(\scA, m)$.  
\end{lemma}

\proof 
We may assume the coordinate $(x_1, \ldots, x_\ell)$ is generic 
so that no coordinate hyperplane $\{x_i=0\}$ is contained in $\scA$. 
From the assumption, $d\alpha_H\wedge x_i\omega$ has no 
poles along $H$, obviously, so does $d\alpha_H\wedge\omega$. 
Hence $\omega\in\Omega^p(\scA, m)$. 
For $D^p(\scA, m)$ the proof is similar. \qed

\begin{lemma}
\label{lem:shift}
Let $(\scA, m)$ be as in Theorem \ref{thm:shift}. 
\begin{eqnarray*}
D^p(\scA, 2k+2\pm m)&\cong& D^p(\scA, 2k\pm m)(-ph), \mbox{ and}\\
\Omega^p(\scA, 2k+2\pm m)&\cong& \Omega^p(\scA, 2k\pm m)(ph).
\end{eqnarray*}
\end{lemma}
\proof 
We only give a proof for $\Omega^p$. The other case is 
immediate from the fact that $D^p$ and $\Omega^p$ are 
dual $S$-modules to each other. 

The case $p=1$ is obvious 
from Theorem \ref{thm:isom} and \ref{thm:der}. 
Put $m'=2k\pm m$. 
Consider the coherent sheaf 
$\scE^p(\scA, m'):=\widetilde{\Omega^p(\scA, m')}$ 
on $\bP^{\ell-1}=\Proj S$ corresponding to the graded 
$S$-module $\Omega^p(\scA, m')$ (\cite{mus-sch}). Recall that 
$\scE^p(\scA, m')$ is known to be a reflexive $\scO$-module, and 
from Lemma \ref{lem:sat} $\Omega^p(\scA, m')$ can be recovered 
from $\scE^p(\scA, m')$ by taking the global section 
$\Gamma_*(\scE^p(\scA, m')):=\bigoplus_{d\in\bZ}
\Gamma(\bP^{\ell-1}, \scE^p(\scA, m')(d))=\Omega^p(\scA, m')$. 
Let $L(\scA)$ be the intersection lattice, and denote 
by $L_k(\scA)$ the set of intersections of codimension $k$. 
For $X\in L_2(\scA)$, denote by 
$\overline{X}\subset\bP^{\ell-1}$ the corresponding flat. 
Consider the open subset 
$$
U=\bP^{\ell-1}\setminus 
\bigcup_{X\in L_2(\scA)}\overline{X}
$$
with the inclusion $i:U\hookrightarrow\bP^{\ell-1}$. 
Since $\scE^p(\scA, m')$ is reflexive, hence normal, 
we have $i_*\scE^p(\scA, m')_U\cong \scE^p(\scA, m')$. 
Furthermore, since $\scE^p(\scA, m')_U$ is locally free on $U$, 
we have 
$$
\scE^p(\scA, m')_U\cong \wedge^p \scE^1(\scA, m')_U. 
$$
Combining these facts, we have 
\begin{eqnarray*}
\scE^p(\scA, m'+2)&=&i_*\scE^p(\scA, m'+2)_U\\
&=&i_*\left(\wedge^p\scE^1(\scA, m'+2)_U \right)\\
&=&i_*\left(\wedge^p\left(\scE^1(\scA, m')_U\otimes \scO(h)_U\right) \right)\\
&=&i_*\left(\scE^p(\scA, m')_U\otimes\scO(ph)_U \right)\\
&=&\scE^p(\scA, m')\otimes\scO(ph). 
\end{eqnarray*}
By taking the global section, we have 
$\Omega^p(\scA, 2k+2\pm m)\cong \Omega^p(\scA, 2k\pm m)(ph)$. 
\qed
\medskip

\noindent
\textit{Proof of Theorem \ref{thm:shift}.} 
Let us prove (2). 
From Theorem \ref{thm:isom} and Lemma \ref{lem:shift}, 
we obtain the isomorphism 
$\Omega^p(\scA, 2k-m)\cong D^p(\scA, m)(pkh)$ of 
graded $S$-modules. Hence their Hilbert series are related by 
the relation 
$$
H(\Omega^p(\scA, 2k-m), q)=H(D^p(\scA, m), q)q^{-pkh}. 
$$
From the definitions of $\phi$ and $\psi$, 
\begin{eqnarray*}
\phi(\scA, 2k-m; t,q)&=&\sum_{p=0}^\ell H( \Omega^p (\scA, 2k-m), q)(t(1-q)-1)^p\\
&=&\sum_{p=0}^\ell H(D^p(\scA, m), q)q^{-pkh}(t(1-q)-1)^p\\
&=&\sum_{p=0}^\ell H(D^p(\scA, m), q)\{q^{-kh}(t(1-q)-1)\}^p, \\
&=&\psi(\scA, m; \frac{q^{-kh}-1}{1-q}-q^{-kh}t,q). 
\end{eqnarray*}
Now we have $\phi(\scA, 2k-m; t, 1)=\psi(\scA, m; kh-t, 1)$ as $q\rightarrow1$ 
and obtain (2). The proof of (1) is similar. \qed

\begin{example}
Suppose $\scA$ is defined by 
$xyz(x+y)(y+z)(x+y+z)$, which is linearly isomorphic to 
the Coxeter arrangement of type $A_3$ and $h=4$. Let 
$m:\scA\rightarrow\{0,1\}$ be defined by $m^{-1}(1)=xyz(x+y+z)$. 
Then $\chi((\scA, m), t)=t^3-4t^2+6t-3$. 
Thus we have from Theorem \ref{thm:shift} that 
\begin{eqnarray*}
\chi ( ( \scA, 2k+m ), t )&=& (t-4k)^3-4(t-4k)^2+6(t-4k)-3\\
\chi ( ( \scA, 2k-m ), t )&=& (t-4k)^3+4(t-4k)^2+6(t-4k)+3. 
\end{eqnarray*}
\end{example}

Theorem \ref{thm:shift} says that for any quasi-constant 
multiplicity $m'$ on a Coxeter arrangement $\scA$ with 
the Coxeter number $h$, the formula 
$$
\chi((\scA, m'+2k+2), t)=\chi((\scA, m'+2k), t-h)
$$
holds. Some computational examples show that 
similar formula holds 
for any multiplicity $m':\scA\rightarrow\bZ_{\geq 0}$, namely, 
supporting the following conjecture.

\begin{conjecture}
Let $\scA$ be a Coxeter arrangement with the Coxeter number $h$. 
Let $m:\scA\rightarrow\bZ_{\geq 0}$ be a multiplicity. 
Then there exists a constant $N=N(\scA, m)$ such that 
$$
\chi((\scA, m+2k+2), t)=\chi((\scA, m+2k), t-h) 
$$
is satisfied for any integer $k>N$. 
\end{conjecture}

\medskip

{\bf Acknowledgment}
T. A. is supported by 21st Century COE program ``Mathematics of Nonlinear 
Structures via Singularities'' Hokkaido University. 
M. Y. is supported by JSPS Postdoctoral Fellowship 
for Research Abroad. 
This paper was begun when T. A. was staying at the 
Adbus Salam ICTP as a short time visitor. 
Both authors are grateful to the Abdus Salam ICTP 
for the hospitality.

\noindent
Takuro Abe\\
Department of Mathematics, \\
Hokkaido University, \\
Sapporo 060--0810, Japan. \\
abetaku@math.sci.hokudai.ac.jp



\bigskip

\noindent
Masahiko Yoshinaga\\
The Abdus Salam ICTP,\\
Strada Costiera 11, 
Trieste 34014, Italy\\
myoshina@ictp.it


\begin{thebibliography}{99}



\bibitem{atw-char}
T. Abe, H. Terao and M. Wakefield, 
The characteristic polynomial of a multiarrangement. 
Adv. in Math. \textbf{215} (2007), 825--838.



\bibitem{ch}
C. Chevalley, Invariants of finite groups generated by reflections. 
Amer. J. Math. \textbf{77} (1955), 778--782.

\bibitem{er-ab}
P. H. Edelman, V. Reiner, 
Free hyperplane arrangements between $A_{n-1}$ and $B_n$. 
Math. Z. {\bf 215} (1994), 347--365. 

\bibitem{er-rhomb}
P. H. Edelman, V. Reiner, 
Free arrangements and rhombic tilings. 
Discrete Comput. Geom.  {\bf 15}  (1996),  no. 3, 307--340. 

\bibitem{mus-sch}
M. Musta\c t\v a and H. Schenck, The module of logarithmic $p$-forms
of a locally free arrangement. J. Algebra \textbf{241} (2001),
699--719. 


\bibitem{ot-arr} P. Orlik and H. Terao, Arrangements of hyperplanes. Grundlehren der Mathematischen Wissenschaften, {\bf 300}. Springer-Verlag, Berlin, 1992. xviii+325 pp.

\bibitem{sai-log}
K. Saito, 
Theory of logarithmic differential forms and logarithmic vector fields.
J. Fac. Sci. Univ. Tokyo Sect. IA Math. {\bf 27} (1980), no. 2, 265--291. 

\bibitem{sai-lin}
K. Saito, 
On a linear structure of the quotient variety by a finite reflexion
group. 
Publ. Res. Inst. Math. Sci. {\bf 29} (1993), no. 4, 535--579. 

\bibitem{sai-unif}
K. Saito, 
Uniformization of the orbifold of a finite 
reflection group. {\it Frobenius manifolds}, 265--320, 
Aspects Math., E36 Vieweg, Wiesbaden, 2004. 


\bibitem{st-stf}
L. Solomon, H. Terao,  A formula for the characteristic polynomial of 
an arrangement.  Adv. in Math.  \textbf{64}  (1987),  no. 3, 305--325.

\bibitem{st-double}
L. Solomon, H. Terao,  
The double Coxeter arrangement. Comm. Math. Helv. {\bf 73} (1998) 237--258. 

\bibitem{sta-super}
R. P. Stanley, Supersolvable lattices. 
Algabra Universalis {\bf 2} (1972), 197--217. 



\bibitem{ter-multi}
H. Terao, 
Multiderivations of Coxeter arrangements.  Invent. Math.  {\bf 148}  (2002),  no. 3, 659--674. 

\bibitem{ter-hodge}
H. Terao, The Hodge filtration and the contact-order filtration 
of derivations of Coxeter arrangements.  Manuscripta Math.  {\bf 118} 
(2005),  no. 1, 1--9. 

\bibitem{ter-bases}
H. Terao, Bases of the contact-order filtration 
of derivations of Coxeter arrangements. 
Proc. A. M. S. {\bf 133} (2005), no. 7, 2029--2034. 





\bibitem{yos-multi}
M. Yoshinaga, 
The primitive derivation and freeness of multi-Coxeter 
arrangements. Proc. Japan Acad. Ser A {\bf 78} (2002), 
no. 7, 116--119. 

\bibitem{yos-char}
M. Yoshinaga, Characterization of a free arrangement and conjecture of Edelman and Reiner.  Invent. Math.  {\bf 157}  (2004),  no. 2, 449--454.


\bibitem{zie-multi}
G. Ziegler, 
Multiarrangements of hyperplanes and their freeness.  Singularities (Iowa City, IA, 1986),  345--359,
Contemp. Math., {\bf 90}, Amer. Math. Soc., Providence, RI, 1989. 




\end{thebibliography}
\end{document}